\documentclass[11pt]{article}
\title {The Bennequin number of $n$-trivial closed $n$-braids is negative}
\author {Oliver T. Dasbach \\ University of California,
    Riverside \\Department of Mathematics\\Riverside, CA 92521 - 0135
\and Xiao-Song Lin \\ University of California, Riverside 
              \\Department of Mathematics\\Riverside, CA 92521 - 0135
} 
\author{Oliver T. Dasbach and Xiao-Song Lin
\\University of California,
    Riverside \\Department of Mathematics\\Riverside, CA 92521 - 0135}
\date{}

\usepackage{amsmath,amssymb,amsthm, graphics} 

\newtheorem{theorem}{Theorem}[section]

\newtheorem{proposition}[theorem]{Proposition}

\newtheorem{corollary}[theorem]{Corollary}

\newtheorem{remark}[theorem]{Remark}

\newcommand{\Z} {\mathbb{Z}}
 
\newcommand{\s}{{\sigma}}

\newcommand{\BeT} {\beta_t}   
\newcommand{\br} {\mathfrak{b}}     
\newcommand{\clbr} {\hat{\br}}      

\newcommand{\Pol}{\Gamma}
\newcommand {\gauss}[2]{{\left \lfloor \frac{#1}{#2} \right \rfloor}}

\begin{document} 
\maketitle
\begin{abstract} 
A famous result of Bennequin states that 
for any braid representative of the unknot the 
Bennequin number is negative. We will extend this result to all
$n$-trivial closed $n$-braids. This is a class of infinitely 
many knots closed under taking mirror images.
Our proof relies on a non-standard parametrization of the
Homfly polynomial.
\end {abstract}

\section{Introduction}
The theory of Vassiliev knot invariants leads to a decreasing filtration
on the set of knots. On the other hand, every complexity measurement 
on knots, like the crossing number or the braid index or the minimal genus, 
gives rise to an increasing filtration. 
Although the interplay between these filtrations in opposite directions
on the set of knots is still fairly unclear, the results of this paper 
show that such an interplay shall have a rich content.

In his seminal paper \cite{Bennequin}, Bennequin proved that
for any representative of the unknot as a closed $n$-braid with
exponent sum $e$, it always holds that $e-n<0$. Bennequin's framework
for the proof of this inequality was the study of contact structures
on the 3-space. A smooth knot, transverse to
the standard contact structure on the 3-space, 
has a well-defined self-linking  number. This self-linking number is
an  invariant of transversal knot types.
Bennequin showed that every transversal knot is transversal isotopic
to a closed braid.  In the case when a closed braid $\hat {\br}$ is 
considered as a transverse knot, this self-linking number 
is equal to the {\it Bennequin number} $$\BeT(\br):=e-n.$$

Bennequin's inequality $e-n<0$ for the unknot also follows from work
of Morton \cite{Morton:Seifert_circles} and independent work of 
Franks and Williams  \cite{FW:Bennequin_estimate}. They showed that
for a closed braid $\hat {\br}$ the Bennequin number is always less than the
lowest degree of $v$ in the Homfly polynomial. Since the Homfly
polynomial of an unknot is identically $1$, the Bennequin number is
negative in this case.

Our aim is to show that Bennequin's inequality has a natural finite type 
generalization. We show that an $n$-trivial closed $n$-braid  
has negative Bennequin number (Theorem~\ref{TH}). 
Since taking the mirror image preserves the
$n$-triviality, we actually have $|e|<n$ for an $n$-trivial closed 
$n$-braid. From the point of view of the theory of Vassiliev knot 
invariants, the reason behind this inequality is that algebraically 
independent Vassiliev invariants become dependent on each other when the
braid index is fixed.
 
One motivation of investigating such an interplay between Vassiliev 
invariants and the increasing
filtration on the set of knots given by the braid index
is best described in terms of complexity theory of
Vassiliev invariants coming from the Homfly polynomial.

As shown by Bar-Natan \cite{BarNatan:PolynomialTime}, for a
knot with a diagram of $c$ crossings the computational complexity for
evaluating a given Vassiliev invariant of order $k$ at this
knot is in $O(c^k)$. Moreover, as it was proved in \cite{JVW} the
evaluation of the Jones polynomial at all but eight points is
$\#P$-hard. Recall that $\#P$ is a class containing $NP$.
Since the Jones polynomial is a certain evaluation of the Homfly
polynomial, this results carries over to the Homfly polynomial. There
are points where the evaluation of the Homfly polynomial is $\# P$-hard.

Assuming that $\#P\neq NP$, this implies that for these Vassiliev
invariants, which come from the Homfly polynomial, there cannot
be an upper bound $k$ such that the computation of each of them
is possible in $O(c^k)$-time.

However, the Homfly polynomial has a definitions in terms of a trace
of a certain representation of the braid groups \cite{Jones2}.
This means that for a fixed braid index, the Homfly polynomial
is computable in polynomial time in the word length of the braids. 

Therefore, the braid index must lead to some restrictions for
Vassiliev invariants involving the word length (exponent sum). 
Some of them are outlined below. Some
others are described in \cite{DLL:QuantumMorphing}.

The first author would like to thank Joan Birman for many discussions
on the topics of this paper.

\section {Parametrizing the Homfly polynomial 
as an honest two variable polynomial}

We start with a certain, non-standard definition of the Homfly
polynomial, which seems to be more useful from the Vassiliev
theory point of view.

\begin {proposition}
Let a polynomial $\Pol = \Pol(\br) \in \Z[\mu, z]$ for a braid $\br \in B_n$ 
be defined by
\begin{eqnarray}
\Pol(\br \sigma_i)-\Pol(\br \sigma_i^{-1})&=&z \Pol(\br) \label{skein}\\
\Pol(\br \sigma_n^{-1})&=& (1-\mu z) \Pol(\br)\\
\Pol(\br \sigma_n)&=& \Pol(\br)\\
\Pol(id_n)&=&\mu^{n-1}
\end{eqnarray} 
and $\Pol(\mu, z)$ is invariant under conjugation and relations in the
braid group.
\end{proposition}

\begin{remark}\rm
For a knot $\clbr$ that is the closure of a braid $\br \in B_n$ 
the relation with the Homfly polynomial $P(v,z)$ 
(see e.g. \cite{Lickorish:Knotpolynomials})
is given by:
\begin{eqnarray}\label{Homfly Relation}
\left . \left( (1- \mu z)^{(e-n+1)/2} \Pol (\br)\left(\mu,
     z\right)\right)\right
\vert _ {\mu=\frac{1-v^2}{z}} = P(\clbr)(v,z).
\end{eqnarray}
Here, $e$ denotes the exponent sum of $\br$ with respect to the
standard generators of the braid group.
\end{remark}

Proofs of properties of the polynomial $\Pol(\mu,z)$ rely on the
following way to compute it: If we have a braid $\br$ which closes to
a knot then we unknot the knot by crossing changes. The term on the
right hand-side of the skein relation (\ref{skein}) corresponds to a 
link of two components. 
For a braid $\br$ which closes to a link with $c$ components, we 
unlink it, i.e. we change crossings in two different components. 
The term on the right-hand side of Equation (\ref{skein}) is then
a link with $c-1$ components. 

\begin{proposition} \label{property1}
The polynomial $\Pol(\mu,z)=\sum_j p_j(\mu)z^j$ has the following properties, 
which are easily checked:
\begin{enumerate}
\item For a split link $L=L_1 \cup
L_2$, i.e. $L_1$ and $L_2$ are unlinked,  
\begin{eqnarray}
\Pol(L)= \mu \Pol(L_1) \Pol(L_2).
\end{eqnarray}

\item The degree in $\mu$ of $\Pol(\br)$ for $\br \in B_n$ is less
than $n$, i.e. each of the $p_j(\mu)$'s is a polynomial of degree less than
$n$.  \label{degree in mu:braid restriction}

\item \label{degree in z} The degree in $z$ of $\Pol(\br)$ is less or equal to the word
length of $\br$.
\item The evaluation $\Pol(\br)(0,z)$ is the Alexander polynomial, in
  its Conway form, of the 
closed braid $\hat {\br}$. 
\item The coefficient of $z^k$ in $\Pol(\br)(\mu,z)$, i.e. $p_k$ is a 
Vassiliev knot invariant of order $k$ for the (framed) closed braid $\hat
{\br}$. This follows in the same way as outlined in
\cite{BL,BarNatan1} for related polynomials. 
\end{enumerate}
\end{proposition}

Furthermore we have:
\begin{proposition} \label{property2} 
Assume the braid $\br \in B_n$ 
closes to a link of $c$ components.
\begin{enumerate}
\item \label{evenodd}
The coefficient 
$p_j(\mu)$ is an odd polynomial if $j+c$ is even, otherwise
$p_j(\mu)$ is an even polynomial.
\item \label{degreezero}
$\Pol(\mu,0)=p_0(\mu)=\mu^{c-1}$. In particular $p_0(\mu)=1$ if
the braid closes to a knot.
\item \label{degreeone:knot} If the braid closes to a knot then
$$p_1(\mu)= \frac {e-n+1}{2} \mu,$$ where $e$ is the exponent sum of the braid
$\br$ with respect to the standard generators of the braid group. 


\item For a braid which closes to a knot, the polynomial $p_j(\mu)$ is
of degree less than or equal to $j$. \label{degree in mu}
\end{enumerate}
\end{proposition}

\begin {proof}
All claims follow by an easy induction. Using the skein relation (\ref{skein})
we express a link as the sum or difference of two links which are
either of shorter word length or are ``less linked''. 
\end{proof}

\section{Closed $n$-braids which are $n$-trivial}

We call a knot $n$-trivial if all Vassiliev invariants
up to degree $n$ vanish on it. The unknot is $n$-trivial for each
$n$. It is not known whether there is another knot with this 
property. The Volume
conjecture (see e.g. \cite{Murakamis:VolumeConjecture}) indeed would 
imply that there is none. Even for the Jones polynomial there is no
example of a non-trivial knot with trivial Jones polynomial known
(see e.g. \cite{DH:DoesThe}).

\begin{theorem}\label{TH}
Let $K$ be a $n$-trivial knot which is given as a closed $n$-braid
$\clbr$. Then the Bennequin number $\beta_t(\br)$ is negative.
\end{theorem}

\begin{proof} 
Let $\br$, $\br \in B_n$ be a braid with exponent sum $e$ satisfying 
our conditions.
The polynomial $\Pol(\br)(\mu,z)$ can be written as
$$\Pol(\mu,z)=\sum p_i(\mu) z^i$$
for some polynomials $p_i(\mu)$ of degree less than $n$ in $\mu$.

Now the coefficients of $z^k$ in
$$(1-\mu z)^{(e-n+1)/2} \Pol(\mu,z)$$
is a (framing independent) Vassiliev invariant, which depends on $\mu$ 
and $e$.

Since, by our condition, all Vassiliev invariants up to degree $n$
vanish, we know that $p_i(\mu)$ must be the coefficient
$c_i(\mu,e)$ of $z^i$ in $(1-\mu z)^{-(e-n+1)/2}$.

We know that $p_i(\mu)$ must have degree less than $n$, but
$c_n(\mu,e)$ has degree $n$. Thus $c_n(\mu,e)$ must be trivial
for all choices of $\mu$, which means, $e$ has to be a root of
$c_n(\mu,e)$.

Hence, by applying the Taylor expansion formula to 
$(1-\mu z)^{-(e-n+1)/2}$, we conclude
that $e$ is one of the numbers $n-1, n-3, \dots, n-2n+1=-n+1$.
\end{proof}

\begin{remark}{\rm
(Unfortunately) the fact that a closed $n$-braid is $n$ trivial
does not imply that it is trivial or even that its Homfly polynomial is
trivial. By the work of Stanford \cite{Stanford} 
if two braids differ by an element in the $(k+1)$-th term of the 
lower central series of the pure braid group
$P_n$ then all Vassiliev invariants up to degree $k$ coincide.
Now, take e.g. the braid group $B_3$ on $3$ strands. The pure braid
group $P_3$ is the direct product of the center of $B_3$ and a
free group of rank $2$, generated by $\s_1^2$ and $\s_2^2$. 
(This is a special case of Newworld's Lemma \cite{DM:Subgroup_separable}.)
By the classification of links given as closed $3$-braids
\cite{BirmanMenasco:3braids} the closure of the 
concatenation of the braid $\s_1 \s_2$ with an element in the $(k+1)$-th
term of the lower central series
is a non-trivial knot. By Stanford's work, it is $k$-trivial, though.
Since the Jones polynomial of a closed $3$-braids is identically $1$, if and 
only if the knot is trivial \cite{Birman}, 
it follows that the Jones polynomial of this link and thus the 
Homfly polynomial is non-trivial. }
\end{remark}

\section{Dimensions of Homfly subspaces restricted to $n$-braids}

We would like to point out one fact which becomes more apparent in the
framed and reparametrized version $\Gamma(\mu,z)$ of the Homfly
polynomial. Vassiliev invariants are dual to the space of knots.
Therefore, one can ask about the dimension of Vassiliev invariants
on subsets of knots. 

We get 
\begin{proposition} 
\begin{enumerate}
\item \label{partone}
The dimension of (framed) Vassiliev knot invariants coming from the
coefficient of the (framed) Homfly polynomial $\Gamma(\mu,z)$ has
dimension
$\gauss k 2 +1$ in degree $k$.

\item
Restricted to closed braids in $B_n$ one gets for the dimension of
Vassiliev invariants in degree $k$:

\begin{enumerate}
\item If $k<n$ then the dimension is $\gauss k 2 +1$.
\item If $k \geq n$ and $n$ is odd and 
$k$ is even then the dimension is $\gauss n 2+1$.
In all other cases it is $\gauss n 2$.
\end{enumerate}
\end{enumerate}
\end{proposition}

\begin{proof}
That the given dimensions form an upper bound, follows from 
Proposition \ref {property1} (\ref{degree in mu:braid restriction})
and Proposition \ref {property2} (\ref{evenodd}) and (\ref{degree in mu}).

It remains to construct the following: Let $k$ be the degree. We will
show that there are braids in $B_{k+1}$ such that the dimension of the
coefficients of $z^k$ is $\gauss k 2 +1$.

For a braid $\br \in B_j$ the polynomial $\Pol(\mu,z)$ satisfies:
\begin{eqnarray}
\Pol(\br \s_j^3)(\mu,z)&=&(1+z \mu + z^2) \Pol(\br)(\mu,z).
\end{eqnarray}

Furthermore, we have for arbitrary $k$:
\begin{eqnarray}
\Pol(\s_1^{-1} \s_2^{-1} \dots \s_k^{-1})(\mu,z)&=&(1-\mu z)^k.
\end{eqnarray}
In particular, the coefficient of $z^k$ is $\pm \mu^k$.

Combining these two formulas, one sees that the coefficient of $z^k$ in
$$\Pol(\s_1^{-1} \s_2^{-1} \dots \s_{k-2}^{-1} \s_{k-1} \s_k^3)(\mu,z)$$
is $\pm \mu^{k-2}$. Replacing successively pairs $\s_{j-1}^{-1}
\s_j^{-1}$ by $\s_{j-1} \s_j^3$, where $j$ has the same parity as $k$, 
yields polynomials with coefficients of $z^k$ equal to $\pm \mu^{j-2}$. 
  
Hence the dimension is equal to the number of $j \leq k$ with the same 
parity as $k$.
\end{proof}

\begin{remark} \rm
Part (\ref{partone}) of the previous proposition should also follow
along the lines of \cite{Meng}.
\end{remark}

\section{Knots whose evaluation of Vassiliev invariants coincide}

Our reparametrization allows us to give an improvement of a theorem
in \cite{KSS:finiteness}. The general question is: Given a quantum
polynomial and a knot $K$ of crossing number $c$, up to which degree
- as a function in $c$ - does one has to know the values of Vassiliev 
invariants at the knot $K$, so that the whole polynomial is already
determined? Since there are only finitely many knots of crossing
number $c$, such a function has to exist. It is not clear, however,
whether this function has some nice form.  In \cite{KSS:finiteness} it
was shown that for the Homfly polynomial such a function is bounded by
a quadratic polynomial in $c$. We will show, that it is in fact
bounded by $c$. Why is this important? It is still unknown whether 
Vassiliev invariants can distinguish knots. The
number of knots of crossing number $c$ grows exponential 
\cite{Welsh:Number_of_knots}. On the other hand, the best lower bound
known for the dimensions of Vassiliev invariants of degree $c$ is
exponential in the square root of $c$ \cite{Dasbach4}.
This might indicate that the space of Vassiliev invariants is
simply not big enough to distinguish knots.

\begin{theorem}
Let the knots $\clbr_1$ and $\clbr_2$ be closures of braids
of length less or equal than $c$, for some number $c$.

If the evaluation of all Vassiliev invariants up to order $c$
coincide on $\clbr_1$ and $\clbr_2$ then their Homfly polynomials are
equal.
\end{theorem}

\begin{proof}
By Proposition \ref{property1} (\ref {degree in z}) we know that
the polynomials $\Gamma (\br_1)$ and $\Gamma (\br_2)$ are polynomials
in $z$ of degree less than or equal to $c$. 

Let $e_1$ ($e_2$, respectively) be the exponent sum of the braid $\br_1$
($\br_2$, respectively) on $n_1$ ($n_2$, respectively) strands. 

We know that the Homfly polynomial $P(\mu,z)$ - with our reparametrization -
is given by a power series in $z$
$$
P(\clbr_1)(\mu,z)= (1-\mu z)^{(e-n+1)/2} \Gamma(\br_1)(\mu,z).
$$

The coefficient of $z^j$ in $P(\clbr_1)(\mu,z)$ is a Vassiliev 
invariant of order $j$.

Since $\Gamma(\clbr_1)(\mu,z)$ is of degree in $z$ less than or equal to
$c$, $\Gamma(\br_1)(\mu,z)$ is determined by all Vassiliev invariants 
up to order $c$. In turn,  $P(\clbr_1)(\mu,z)$ itself is determined
by all terms in $z^j$ of degree less than or equal to $c$.

The same argument works with $\br_2$ and we are done.
\end {proof}

\begin{remark} \rm
For the definition of $\Gamma$ we do not really need the knot given
as a closed braid. We can work with arbitary knot diagrams of framed 
knots instead. Thus, in the last theorem we can replace the word length 
by the crossing number. As a corollary we get:
\end{remark}

\begin{corollary}
If all Vassiliev invariants up to degree $c$ vanish on 
a knot $K$ of crossing number $c$ then the knot has trivial Homfly
polynomial.
\end{corollary}

\providecommand{\bysame}{\leavevmode\hbox to3em{\hrulefill}\thinspace}


\begin{thebibliography}{JVW90}

\bibitem[Ben83]{Bennequin}
D.~Bennequin, \emph{Entrelacement et \'equations de {P}faff}, Ast\'erisque
  \textbf{107 - 108} (1983), 83 -- 161 (French).

\bibitem[Bir85]{Birman}
J.~S. Birman, \emph{{ On the Jones Polynomial of Closed 3-braids}}, Invent.
  Math. \textbf{81} (1985), 287--294.

\bibitem[BL93]{BL}
J.~S. Birman and X.-S. Lin, \emph{Knot polynomials and {V}assiliev's
  invariants}, Invent. Math. \textbf{111} (1993), no.~2, 225--270.

\bibitem[BM93]{BirmanMenasco:3braids}
J.~S. Birman and W.W. Menasco, \emph{Studying links via closed braids {III}.
  {C}lassifying links which are closed 3-braids}, Pacific J. Math. \textbf{161}
  (1993), no.~1, 25--113.

\bibitem[BN95a]{BarNatan:PolynomialTime}
D.~Bar-Natan, \emph{{ Polynomial Invariants are Polynomial}}, Math. Research
  Letters \textbf{2} (1995), 239--246.

\bibitem[BN95b]{BarNatan1}
D.~Bar-Natan, \emph{On the {V}assiliev knot invariants}, Topology \textbf{34}
  (1995), 423--472.

\bibitem[Das00]{Dasbach4}
O.T. Dasbach, \emph{On the combinatorial structure of primitive {V}assiliev
  invariants {III} - {A} lower bound}, Comm. Contempor. Math. \textbf{4}
  (2000), to appear, available as: math.GT/9806086.

\bibitem[DH97]{DH:DoesThe}
O.T. Dasbach and S.~Hougardy, \emph{{Does the Jones Polynomial Detect
  Unknottedness?}}, Experimental Math. \textbf{6} (1997), no.~1, 51 -- 56.

\bibitem[DLL00]{DLL:QuantumMorphing}
O.T. Dasbach, T.D. Le, and X.-S. Lin, \emph{Quantum morphing and the {J}ones
  polynomial}, preprint, 2000.

\bibitem[DM99]{DM:Subgroup_separable}
O.T. Dasbach and B.~Mangum, \emph{The automorphism group of a free group is not
  subgroup separable}, to appear in: Knots, Braids, and Mapping Class Groups
  Conference Proceedings, Internatinoal Press.

\bibitem[FW87]{FW:Bennequin_estimate}
J.~Franks and R.~F. Williams, \emph{Braids and the {J}ones polynomial}, Trans.
  Am. Math. Soc. \textbf{303} (1987), 97--108.

\bibitem[Jon87]{Jones2}
V.~F.~R. Jones, \emph{{H}ecke algebra representations of braid groups and link
  polynomials}, Ann. of Math. \textbf{126} (1987), 335--388.

\bibitem[JVW90]{JVW}
F.~Jaeger, D.~L. Vertigan, and D.~J.~A. Welsh, \emph{On the computational
  complexity of the {J}ones and {T}utte polynomials}, Math. Proc. Cambridge
  Philos. Soc. \textbf{108} (1990), no.~1, 35--53.

\bibitem[KSS97]{KSS:finiteness}
L.H. Kauffman, M.~Saito, and S.F. Sawin, \emph{On finiteness of certain
  {V}assiliev invariants}, J. Knot Theory Ramifications \textbf{6} (1997),
  no.~2, 291--297.

\bibitem[Lic88]{Lickorish:Knotpolynomials}
W.B.R. Lickorish, \emph{Polynomials for links}, Bull. Lond. Math. Soc.
  \textbf{20} (1988), 558--588.

\bibitem[Men97]{Meng}
G.~Meng, \emph{Bracket models for weight systems and the universal {V}assiliev
  invariants}, Topology Appl. \textbf{76} (1997), no.~1, 47--60.

\bibitem[MM99]{Murakamis:VolumeConjecture}
H.~Murakami and J.~Murakami, \emph{The colored {J}ones polynomials and the
  simplicial volume of a knot}, preprint, available as: math.GT/9905075, 1999.

\bibitem[Mor86]{Morton:Seifert_circles}
H.R. Morton, \emph{Seifert circles and knot polynomials}, Math. Proc. Camb.
  Phil. Soc. \textbf{99} (1986), 107--109.

\bibitem[Sta96]{Stanford}
T.~Stanford, \emph{Braid commutators and {V}assiliev invariants}, Pacific J.
  Math. \textbf{174} (1996), no.~1, 269--276.

\bibitem[Wel92]{Welsh:Number_of_knots}
D.J.A. Welsh, \emph{On the number of knots and links}, Sets, graphs and numbers
  (Budapest, 1991), North-Holland, Amsterdam, 1992, pp.~713--718.

\end{thebibliography}
\end{document}